\documentclass[12pt,a4paper]{article}

\usepackage{color,srcltx,latexsym,amssymb,amsmath,amsthm}

\newcommand{\pr}{{\bf P}}                     
\newcommand{\exn}{{\bf E\,}}                    
\newcommand{\ind}{{\bf 1}}                     
\newcommand{\ep}{\varepsilon}            
\newcommand{\ept}{{\varepsilon_T}}            
\newcommand{\la}{\lambda}            
\newcommand{\al}{\alpha}                 
\newcommand{\de}{\delta}                 
\newcommand{\ze}{\zeta}                 

\newcommand{\Bom}{B_M}                 

\newcommand{\A}{{A}}      
\newcommand{\Ao}{{A_0}}      

\newcommand{\Aw}{{\widehat A}}        
\newcommand{\Dw}{{\widehat D}}        

\newcommand{\lal}{\mathcal L}                 

\newcommand{\R}{\mathbb{R}}             

\newcommand{\oo}{\overline{o}}
\newcommand{\ui}{\underline\infty}

\newcommand{\cl}[1]{\operatorname{cl}(#1)}
\newcommand{\intr}[1]{\operatorname{int}(#1)}

\newcommand{\doma}{\mbox{\rm dom}\, \Ao}
\newcommand{\domd}{\mbox{\rm dom}\, D}

\newcommand{\grad}{\mathop{\rm grad}}

\newtheorem{Theorem}{Theorem}

\newtheorem{Def}{Defintion}
\newtheorem{Rem}{Remark}
\newtheorem{Example}{Example}

\newcommand\Proof{{\sl Proof}}

\title{On necessary and sufficient conditions for the local large deviation principle}

\author{Konstantin Borovkov.\footnote{School of Mathematics and Statistics, The University of Melbourne, Parkville VIC 3010, Australia; e-mail: borovkov@unimelb.edu.au.
}}

\date{}

\begin{document}

\maketitle

\begin{abstract}
One says that the local large deviation principle (LLDP) is satisfied for a   family of random vectors $\{\zeta_T\}_{T\ge 0}$ in~$\R^d,$ $d\ge 1,$ if there exists a function $D:\R^d\to [0,\infty],$ $D\not \equiv \infty,$ such that, for any $\al\in \R^d$,
\[
\lim_{T\to\infty} T^{-1}\ln \pr (|\zeta_T -\alpha |<\ep_T) = - D(\al)
\]
for $\ep_T\to 0$ slowly enough. In this paper, we establish necessary and sufficient conditions for the LLDP that are very close to each other. Namely, if the LLDP is satisfied then, for $M_T\to\infty$ slowly enough as $T\to\infty$, there exists the limit
\[
\A(\mu):= \lim_{T\to\infty}T^{-1}\ln \exn (e^{T\langle \mu, \zeta_T\rangle}; |\zeta_T|\le M_T)\in (-\infty, \infty] , \quad \mu\in \R^d,
\]
which is equal to the Legendre--Fenchel transform $\mathcal L_D$ of the rate function~$D$.
Conversely, if the above limit~$\A(\cdot )$ exists and is an essentially smooth function, then the LLDP is satisfied with the rate function~$D$ equal to~$\mathcal L_A.$ This ``relaxed version'' of the G\"artner--Ellis theorem's main condition does not involve the restrictive integrability assumptions from the latter and is most adequate to the nature of the local large deviation problem.

\bigskip \noindent
{\em Key words and phrases}: large deviation principle, local large deviation principle, Lagrange--Fenchel transform, G\"artner--Ellis theorem, Laplace--Varadhan lemma.

\bigskip  \noindent
{\em AMS 2020 Subject Classification}: Primary 60F10; secondary 26B25.
\end{abstract}

\section{Introduction and the main result}

Systemic studies of the large deviation (LD) probabilities, in the first place for the classical scheme of summation of independent identically distributed  (IID)
random variables, goes back to  the seminal paper~\cite{Cr38}. Since then, a vast literature devoted to the LD theory, both in the classical and more general settings,  has emerged. Surveys of the history of the development of that theory and its fundamental results can be found in~\cite{BoMo92, Bo20a, BoMo11, DeZe09, DeSt89, Pu01}.

A standard framework for characterising the asymptotic behaviour of the probabilities of rare events is provided by the  {\em large deviation principle\/} (LDP;  this term was apparently first used in~\cite{Va66}). We will state it for a general family $\{\ze_T\}_{T\ge 0}$ of random vectors in $\R^d,$ $d\ge 1,$ where $T$ can be a discrete ($T\in \{1,2,\cdots\}$) or continuous ($T\in  [0,\infty)$) parameter.

\begin{Def}
{\em
One says that a family  $\{\ze_T\}_{T\ge 0}$ of random vectors in~$\R^d$ satisfies the LDP if there exists a lower  semicontinuous function   $D:\R^d\to [0,\infty]$,  $D\not\equiv \infty$, such that, for any Borel set $B\subset \R^d,$
\begin{align}
\liminf_{T\to\infty}\frac1{T}\ln\pr (\zeta_T\in B) & \ge - D(\intr{B}),
\notag\\
  \limsup_{T\to\infty}\frac1{T}\ln\pr (\zeta_T\in B) & \le - D(\cl{B}),
  \label{upper_LDP}
\end{align}
where
\[
D(B):=\inf_{\al \in B}D(\al ),
\]
and we denoted by $\intr{B}$  the interior of  $B$ and by $\cl{B}$ its closure.
}
\end{Def}

The function  $D$  is referred to as the {\em  rate function\/} for $\{\ze_T\}_{T\ge 0}$. A rate function~$D$ is called ``good'' if
\begin{align}
\label{GoodRF}
\{\al\in \R: D(\al)\le v\} \quad \mbox{is compact for any $v>0$}.
\end{align}

In the classical discrete time random walk LD setup, one deals with the family 
$\{\ze_T\}$  of sample means
\begin{align}
\label{Classic}
\ze_T:=\frac1T \sum_{j=1}^T \xi_j, \quad T=1,2,\ldots,
\end{align}
of  
IID random vectors   $\xi_j\in \R^d,$ $j\ge 1$, that satisfy the moment Cram\'er condition: denoting by $\langle \cdot, \cdot \rangle$ the Euclidean scalar product in~$\R^d$, one   assumes that
\[
\psi (\mu):=\ln \exn e^{\langle \mu ,\xi_1\rangle }<\infty
\]
for $\mu$ from a non-empty neigbourhood of the origin. In this case, the celebrated Cram\'er theorem (see e.g.\ Theorem~2.2.30 in~\cite{DeZe09}) states that   $\{\ze_T\}$ satisfies the LDP with the rate function given by the Legendre--Fenchel transform of the (convex) cumulant generating function~$\psi$:
\[
D(\al )= \lal_\psi (\al):= \sup_{\mu\in \R^d} (\langle \al ,\mu\rangle- \psi (\mu))\in [0,\infty], \quad \alpha \in \R^d.
\]
It is well-known that the transform   $\lal_{\bullet}$  maps the class~$\mathcal C_d$ of convex functions on~$\R^d$ into itself and, moreover, on~$\mathcal C_d$ the mapping~$\lal_{\bullet}$ is inverse to itself   (see e.g.\ \S26 in~\cite{Ro70}). Therefore, in this case, $D$ is also convex and one has   $\psi (\mu) =\lal_D (\mu)$, $\mu\in \R^d$.

Note that, for a family $\{\ze_T\}$ of random vectors of a general nature,  it is possible to have the LDP satisfied  with a non-convex rate function (see e.g.\ p.~75 in~\cite{DeZe09}, and also our Examples~\ref{Ex_2} and~\ref{Ex_3} below).

The probabilistic meaning of the rate function is best seen from the so-called {\em local\/} LDP (LLDP) introduced in~\cite{BoMo10, BoMo11}. For $\al\in \R^d,$ set  $|\al|:=\langle \al, \al \rangle^{1/2}$ and, for $\ep >0,$ denote by
\[
(\al)_\ep:=\{\beta\in \R^d: |\al -\beta|<\ep\}
\]
the  $\ep$-neighbourhood of~$\al.$

\begin{Def}\label{def_lpbu}
{\em
One says that a family  $\{\zeta_T\}_{T\ge 0}$ of random vectors in $\R^d$ satisfies the LLDP if there exists a function $D:\R^d\to [0,\infty]$ such that $D\not \equiv \infty$ and, for any $\al \in \R^d,$ one has
\begin{align}
\lim_{\ep\searrow 0}\liminf_{T\to\infty}\frac1{T}\ln\pr (\zeta_T\in (\al)_\ep)
 & =- D( \al ),
\label{loc_ld_inf}
  \\
\lim_{\ep\searrow 0}\limsup_{T\to\infty}\frac1{T}\ln\pr (\zeta_T\in (\al)_\ep)
 & = - D( \al ).
\label{loc_ld_sup}
\end{align}
}
\end{Def}

It was shown in~\cite{BoMo11} that if the LLDP is satisfied then the function~$D$ in Definition~\ref{def_lpbu} is automatically lower semicontinuous.

\begin{Rem}
\rm
Note that the LLDP can be satisfied in situations where all the probability mass but its exponentially small portion escapes to infinity ``at a higher speed''  (so that $\inf_\al D(\al) >0$), which is not possible when the ``usual'' LDP is satisfied. See Theorem~\ref{equiv}(ii) and Example~\ref{Ex_2}(ii) below.
\end{Rem}

To give a convenient equivalent definition of the LLDP, we need to introduce the following concept.

\begin{Def}\label{def_slow}
{\rm
We will say that a statement~$\mathcal S$, whose formulation involves a function $\ep=\ept >0 $, {\em is satisfied for $\ept $ that tends  to zero slowly enough  as\/} $T\to\infty $ (and will denote this property of the function~$\ep$ by writing  ``$\ep =\oo (1)$ as $T\to\infty $'') if
\begin{enumerate}
\item[(i)] there exists a positive function  $\widetilde\ep = \widetilde\ep_T \to 0$ as $T\to\infty  $ such that   $\mathcal S$ holds true with $\ep=\widetilde\ep,$ and
\item[(ii)] the statement~$\mathcal S$ holds true for all functions
$\ept \to 0 $ such that
$\ept \ge   c \widetilde\ep_T $ for some constant  $c>0.$

\end{enumerate}
Likewise, we will say that  a statement~$\mathcal S$, whose formulation involves a function $M=M_T >0, $ {\em is satisfied for $M_T$ that tends  to infinity slowly enough  as\/} $T\to\infty $ (and will denote this property of the function~$M$ by writing  ``$M \to \ui$ as $T\to\infty $'') if~$\mathcal S$ holds when $1/M=\oo (1)$ as $T\to\infty $.
}
\end{Def}

\begin{Rem}
\em
Suppose  $\mathcal S$ is satisfied for $\ep =\oo (1)$ and let $\widetilde\ep$ be the function from~(i) in Definition~\ref{def_slow}. Then~$\mathcal S$ does not need to hold for  a function~$\ep=\widehat\ep$, such that the relation $\widehat\ep_T \ge   c \widetilde\ep_T $  is not true for any $c >0$  (e.g.,  $\ept \ll  \widetilde \ep_T  $ as $T\to\infty $). For instance, let~$\mathcal S$ denote the statement that ``$\ep T\to \infty$ as $T\to\infty$''. It clearly holds true for any fixed $\ep >0$. Furthermore, it is also true for $\ep = \ept=\oo (1)$ (one can take, say, $  \widetilde \ep_T= T^{-1/2}$). However, for $\ept = T^{-1}\ll \widetilde \ep_T$,  the statement~$\mathcal S$ is not true anymore.
\end{Rem}

It was shown in~\cite{BoMo11} that   the LLDP is equivalent to the following: if
$ \ept =\oo(1)$ as $T\to\infty,$ then
\begin{align}
\label{lpbu}
\lim_{T\to\infty}\frac1{T}\ln\pr (\zeta_T\in (\al)_\ept) = - D( \al ), \quad
\al \in \R^d.
\end{align}
Moreover, under the assumption that the respective rate function~$D$ is good (see~\eqref{GoodRF}), Theorem~1.1 in~\cite{BoMo11} implies the following assertion.

\begin{Theorem}
\label{equiv}
{\rm (i)}~In the general case, if  the LDP is satisfied for a family of random vectors then the LLDP with the same rate function~$D$ is also satisfied for that family.

{\rm (ii)}~The converse implication holds true provided that the following condition is met{\rm :}  for any $N>0$ there exists a $v>0$ such that
\begin{align}
\label{About_tail}
\limsup_{T\to\infty} \frac1T \ln \pr (|\zeta_T|> v) \le -N.
\end{align}
\end{Theorem}

The last condition is often referred to as {\em exponential tightness}.
Thus, the LDP and LLDP are equivalent provided that both conditions~\eqref{GoodRF} and~\eqref{About_tail} are satisfied, i.e.\ when the rate function is good and the family $\{\zeta_T\}$ is exponentially tight. When~\eqref{About_tail} is not met, the probability mass ``escapes to infinity at a higher speed'', cf.\ Example~\ref{Ex_2}(ii) below.

In the case of a general family $\{\ze_T\}$ of random vectors, sufficient conditions for the LDP in terms of the behavior of the moment generating function  were established in the celebrated G\"artner--Ellis theorem  (\cite{Ga77,El84}; see also Section~2.3 in~\cite{DeZe09}) to be stated next.
To clarify the assumptions of the theorem, we first recall that a function $f\in \mathcal C_d$ with effective  domain
\[
\mbox{dom}\, f :=  \{\mu\in \R^d: f(\mu ) <\infty\}
\]
is said to be {\em essentially smooth\/} if:

\begin{enumerate}
\item[(i)]   $\intr{\mbox{\rm dom}\, f }\neq \varnothing$,

\item[(ii)]  $f$ is differentiable in $\intr{\mbox{\rm dom}\, f}$,
\smallskip

\item[(iii)] $\lim_{k\to\infty} |\grad f (\mu_k)|=\infty$ for any sequence  $\{\mu_k\}_{k\ge 1}\subset \intr{\mbox{\rm dom}\, f}$ such that there exists the limit $\lim_{k\to\infty} \mu_k =\mu \in \partial\, \mbox{\rm dom}\, f$
\end{enumerate}
(see e.g.\ \S26 in~\cite{Ro70}  or Definition~2.3.5 in~\cite{DeZe09}; note that if  $\mbox{\rm dom}\, f$ is unbounded in some direction, the ``infinitely remote point'' in that direction is not considered to be part of $\partial\, \mbox{\rm dom}\, f$).

\begin{Theorem}
\label{G-E} {\rm (G\"artner--Ellis)}
Assume that, for any $\mu\in \R^d$, there exists the limit
\begin{align}
\label{base_function}
\Ao (\mu):=\lim_{T\to\infty}\frac1T \ln \exn e^{T \langle \mu, \ze_T\rangle } \in (-\infty, \infty]
\end{align}
which  is essentially smooth and, moreover,  $0 \in \intr{\doma}$. Then the LDP is satisfied for $\{\ze_T\}$ with the rate function $D=\mathcal L_\Ao$ for which condition~\eqref{GoodRF} is met.
\end{Theorem}

Following~\cite{Bo20}, we will call the limit $\Ao$ in~\eqref{base_function} the {\em fundamental function} (FF) for $\{\ze_T\}$ and say that it exists if  $  \intr{\doma}\neq \varnothing$.  Note that, since~$\Ao \in \mathcal C_d$ (due to its  being the point-wise limit of convex functions), one has $\Ao=\mathcal L_D$.

There is also an assertion that is converse to the  G\"artner--Ellis theorem and goes back to~\cite{Va66}. It holds under additional integrability assumptions and is often referred to as the Varadhan (or Laplace--Varadhan) integral lemma, see  e.g.\ Theorem~4.3.1 in~\cite{DeZe09}. Namely, if the LDP is satisfied for   $\{\ze_T\}$ with a good rate function~$D$ and, for some $\mu\in \R^d$, one has
\begin{align}
\label{Vara}
\lim_{M\to\infty} \limsup_{T\to\infty} \frac1T
\ln \exn (e^{T \langle \mu, \ze_T\rangle} ;   \langle \mu, \ze_T\rangle\ge  M ) =-\infty,
\end{align}
then, for this $\mu$, there exists the limit~\eqref{base_function} that is equal to~$\mathcal L_D (\mu)$. Moreover, if one replaces $ \langle \mu, \ze_T\rangle$ in~\eqref{Vara}  with  $ \phi (\zeta_T)$, where  $\phi:\R^d \to\R$ is an arbitrary continuous function, then there will also exist the limit
\[
\lim_{T\to\infty} \frac1T \ln \exn e^{T \phi (\zeta_T)  } =\sup_{\al}(\phi (\al)- D(\al)).
\]

Furthermore, for compound renewal processes $\{Z_T\}_{T\ge 0}$, it was established in~\cite{Bo20} that (i)~under suitable assumptions,   the existence of the FF~\eqref{base_function}   for $\ze_T:=T^{-1}Z_T$ implies the LDP for the family  $\{\zeta_T\}$ and, conversely, (ii)~the LDP for $\{\zeta_T\}$ implies the existence of the FF.

Inspection of the above results shows that, for a general family of random vectors
$\{\zeta_T\}$, the existence of the FF and the validity of the LDP are almost equivalent, with the caveat that this relationship requires strong integrability assumptions in both directions (for the direct assertion it is ``built   into''~\eqref{base_function}, while for the converse assertion one requires~\eqref{Vara}). Likewise, the ``transition'' from LLDP to LDP (Theorem~\ref{equiv}(ii)) requires exponential tightness, ensuring that the probabilities of~$\zeta_T$ hitting ``really remote regions'' are negligibly small (condition~\eqref{About_tail}).

However, if one is only after the LLDP or just  the asymptotic behaviour of the probability that $\zeta_T$ hits a given bounded set, then the behaviour of the distribution of $\zeta_T$ in ``really remote regions'' should   be irrelevant. This suggests that the existence of the FF per se is an excessive condition for the  LLDP, as it requires that $\exn e^{T \langle \mu, \ze_T\rangle }<\infty$ for~$\mu$ in an open set. Likewise, condition~\eqref{Vara}, which complements the LDP in the converse implication, basically serves to guarantee  the convergence of the integrals producing in the limit the FF and  is therefore also expected to be excessive for existence of an adequate alternative to the FF in the ``local setup".

The following simple example illustrates the above observations. It shows that the LLDP (and even the LDP) can be satisfied in situations where  $\exn e^{T \langle \mu, \ze_T\rangle } =\infty,$ $\mu \neq 0.$

\begin{Example}\label{Ex_1}
\rm
Let  $\overline{\Phi}(x):=(2\pi )^{-1/2} \int_x^\infty  e^{-y^2/2}dy,$ $x\in \R,$ denote the ``right tail'' of the standard normal distribution, and let $\ze_T$ be a symmetric random variable such that
\begin{align}
\pr (  \ze_T >z)=\left\{
\begin{array}{ll}
\overline{\Phi}(zT^{ 1/2}), & z\in (0, T^{1/2}),
\\
c_T z^{-1},  & z\ge T^{1/2},\quad c_T:= T^{1/2} \overline{\Phi}( T ).
\end{array}
\right.
\label{kontra}
\end{align}
Using   Mills' ratio, we see that, for any fixed $\al >0 , \ep\in (0,\al)$,
for $T> (\al+ \ep)^{1/2}$ one has
\begin{align*}
\pr ( \ze_T \in   (
\al)_\ep) & = \overline{\Phi}((\al -\ep) T^{1/2})-\overline{\Phi}( (\al+ \ep)T^{1/2})
\\
& \sim \overline{\Phi}((\al -\ep) T^{1/2})
 \sim 
 \frac{\exp\{-T\frac{\al^2}2(1 - 2\frac{\ep}{\al} +\frac{\ep^2}{\al^{2}})\}}{(\al -\ep) (2\pi T)^{ 1/2} }  \quad \mbox{as}\ T\to\infty.
\end{align*}
Clearly, the above relations will still hold for $\ep =\ep_T =o(1)$ such that $\ep\gg T^{-1}.$ Hence, for $\{\ze_T\}$, the LLPD is satisfied with the rate function  $D(\al)= \al^2/2,$ $\al \in \R$.

Furthermore, for a fixed $v>0$,  one has from~\eqref{kontra} that, for $T>v^2$,
\[
 \pr ( | \ze_T|> v) = 2 \pr ( \ze_T >v)
  = 2\overline{\Phi}(v T^{ 1/2}) \sim \frac{2^{1/2}e^{-Tv^2/2}}{(\pi T)^{1/2}v}, \quad T\to\infty,
\]
so that  $\lim_{T\to\infty}T^{-1}\ln \pr ( | \ze_T|> v) =  -v^2 /2. $ Therefore, condition~\eqref{About_tail} is met and, by Theorem~\ref{equiv}(ii),  the LDP is also satisfied for $\{\ze_T\}$. At the same time, it is obvious that $\exn e^{T \mu \ze_T} =\infty $ for all  $\mu \neq 0,$ so that the FF does not exist in this example.
\end{Example}

The above suggests considering the following ``relaxed version'' of the G\"artner--Ellis theorem's condition that assumed the existence of the FF.

For $M>0,$ denote by
\[
 \Bom:=\{\al\in \R^d: |\al|\le M\}
\]
the closed ball of radius $M $  centered at the origin and set
\begin{equation}
\label{ETM}
E_{T,M} (\mu):= \exn ( e^{T \langle \mu ,\ze_T \rangle };  \ze_T\in \Bom),
 \quad \mu\in \R^d.
\end{equation}

\begin{Def}\label{def_wsff}
\rm  We will say that there exists a {\it weak sense fundamental function} (WSFF)  $\A(\mu)$  for   family $\{\ze_T\}$ if, for  $M=M_T\to \ui$ as $T\to\infty$,
there exists the limit
\begin{equation}
\label{2}
\lim_{T\to \infty}\frac{1}{T} \ln E_{T,M_T} (\mu) =
\A(\mu) \in (-\infty, \infty],
\quad \mu\in \R^d,
\end{equation}
which is a  convex function with $ \intr{\mbox{dom}\, A }\neq \varnothing$.
\end{Def}

The convexity of $\A$ is automatic, as it is a point-wise limit of convex functions, so that the key requirement here is that the interior of the effective domain of $\A$ is non-empty.

\begin{Rem}
\label{rem_SBF}
\rm
It is not hard to verify that we will get an equivalent definition if, instead of~$\Bom$, we use  in~\eqref{ETM}  sets  of the form     $MB:= \{M\al: \al\in B\} $, where  $B\subset \R^d$ is an arbitrary fixed bounded Borel set such that  $0\in  \intr{B}$.
\end{Rem}

\begin{Rem}
\label{rem_A=A}
\rm
It is obvious  that, although the FF does not exist in Example~\ref{Ex_1}, the WSFF does: for the $\ze_T$'s   with distributions~\eqref{kontra}, relation~\eqref{2} holds true with $\A(\mu) =\mu^2/2,$ $\mu\in\R,$ for any function $M=M_T\to \infty $ as $T\to\infty$ such that $M\le T^{1/2}$. Therefore it holds for $M\to \ui$ as well: in Definition~\ref{def_slow}, for $\ep :=1/M$   it suffices to take, say, $ \widetilde\ep_T =T^{-1/3}.$

Now if FF $\Ao$ does exist for a family $\{\ze_T\}$, then the WSFF $\A$ also exists and
\begin{equation}
\label{A=A}
\A (\mu) = \Ao (\mu ), \quad \mu \in  \intr{\doma}.
\end{equation}
To demonstrate this, assume for simplicity that $d=1$ (the treatment of the case $d>1$ is basically the same). If $\mu \in  \intr{\doma}$ then there exists a $\delta >0$ such that $\mu\pm\delta \in \doma$.  Using~\eqref{base_function} with $\mu \pm \delta$ instead of~$\mu$, we obtain that, for $M_T\to \ui,$ one has
\begin{align*}
\exn ( e^{T\mu \zeta_T}; |\zeta_T | > M_T)
 & =
 \exn ( e^{T\mu \zeta_T};  \zeta_T <- M_T) + \exn ( e^{T\mu \zeta_T}; \zeta_T   > M_T)
 \\
 & \le e^{-T\delta M_T} \big( \exn e^{T(\mu-\delta) \zeta_T}+ \exn e^{T(\mu+\delta) \zeta_T}\big)
 \\
 & =
  e^{-T\delta M_T} \big[  e^{T (\Ao (\mu -\delta) + o(1))}+ e^{T (\Ao (\mu +\delta) + o(1))}\big]
  =  e^{-T\delta M_T(1+o(1))}
\end{align*}
since $|\Ao (\mu \pm \delta) |<\infty.$ Therefore, again appealing to~\eqref{base_function}, we get
\begin{align*}
 \exn ( e^{T\mu \zeta_T}; |\zeta_T | \le M_T)
 &  =  \exn e^{T\mu \zeta_T} - \exn ( e^{T\mu \zeta_T}; |\zeta_T | > M_T)
  \\
 & =
   \exn e^{T\mu \zeta_T}
   \big[1 + O (\exp\{-T [ \delta M_T(1+o(1)) + \Ao (\mu )]\}) \big]
   \\
 &  = \exn e^{T\mu \zeta_T} (1 + o(1)),
\end{align*}
which establishes~\eqref{A=A}.
\end{Rem}

\begin{Rem}
\rm
Akin to the relationship between~\eqref{loc_ld_inf} and \eqref{loc_ld_sup}, on the one hand,  and~\eqref{lpbu} on the other, Definition~\ref{def_wsff} admits an alternative equivalent formulation. Namely, it is not hard to show that the WSFF~$A$ exists iff, for any $\mu\in \R^d$, one has
\begin{equation}
\label{2_eq}
\lim_{N\to\infty} \liminf_{T\to\infty}\frac{1}{T} \ln E_{T,N} (\mu) =
\lim_{N\to\infty} \limsup_{T\to\infty}\frac{1}{T} \ln E_{T,N} (\mu),
\end{equation}
the common value being equal to $\A(\mu)\in (-\infty, \infty].$ We will show that  \eqref{2_eq} implies~\eqref{2}, the converse implication being equally easy to demonstrate.

Suppose that~\eqref{2_eq} holds true and denote the common limit by~$\widetilde A(\mu)$.
It is clear that, for any $M_T\to \infty$ as $T\to\infty$ and fixed $N>0,$ one has $   E_{T,N}  (\mu)\le E_{T,M_T} (\mu)$ for all large enough~$T$. Therefore
\begin{equation}
\label{2_eq0}
\widetilde A(\mu) =
\lim_{N\to\infty }
\liminf_{T\to\infty} \frac{1}{T} \ln E_{T,N}  (\mu)
\le
\liminf_{T\to\infty} \frac{1}{T} \ln E_{T,M_T} (\mu)
.
\end{equation}
It remains to show that, for $M_T\to\ui,$
\begin{equation*}
\limsup_{T\to\infty} \frac{1}{T} \ln E_{T,M_T} (\mu)
\le
\lim_{N\to\infty }
\limsup_{T\to\infty} \frac{1}{T} \ln E_{T,N}  (\mu).
\end{equation*}
To this end, note that, as the right-hand side of~\eqref{2_eq} equals~$\widetilde \A (\mu)$, there exists a sequence $N_k\nearrow \infty$ as $k\to\infty$ such that
\[
\limsup_{T\to\infty}\frac{1}{T} \ln E_{T,N_k} (\mu) \le \widetilde  A (\mu) +k^{-1}/2.
\]
Therefore there exist a sequence $T_k\nearrow \infty$ as $k\to\infty$  such that
\[
\frac{1}{T} \ln E_{T,N_k} (\mu) \le \widetilde A (\mu) +k^{-1} \quad \mbox{for } T\ge T_k.
\]
Now, putting $\overline M_T:= \sum_{k\ge 1} N_k \ind (T\in [T_k, T_{k+1}))$, we obtain a function $\overline M_T\to \infty$ with the property that 
\[
\limsup_{T\to \infty}\frac{1}{T} \ln E_{T,\overline M_T} (\mu)
\le
\widetilde A(\mu).
\]
Together with~\eqref{2_eq0} this means that~\eqref{2} with~$M_T$ replaced by~$\overline M_T$ holds true with $A(\mu)= \widetilde A (\mu)$.
Now setting, say, $\widetilde M_T:= \overline M_T^{1/2}$, it is clear from the above argument that~\eqref{2} still holds for any function  $M_T\le c\widetilde M_T$ for some $c>0$. This means that~\eqref{2} is met for $M_T\to\ui,$ as asserted.
\end{Rem}

It turns out that, modulo natural regularity conditions on~$\A$, the existence of the WSFF is {\em necessary and sufficient\/} for the LLDP to be satisfied. In this assertion, neither the direct nor the converse implications require any integrability conditions, as it should be. Indeed, in both cases, the behaviour of the distributions of the $\zeta_T$'s in ``really remote regions'' (outside the ball $B_M $) is irrelevant for the desired claims.

The following theorem is the main result of this paper.

\begin{Theorem}
\label{t1}
  {\rm (i)} If, for a family  $\{\ze_T\}$, LLDP is satisfied with a rate function $D$ then there exists the WSFF $\A={\mathcal L}_D.$

    {\rm (ii)} If, for a family  $\{\ze_T\}$,  there exists an essentially smooth WSFF $\A$ then the LLDP  with the convex rate function~$D ={\mathcal L}_\A $ is satisfied for that family.
\end{Theorem}

\begin{Rem}
\rm
Note that, unlike the conditions of Theorem~\ref{G-E},   those  of Theorem~\ref{t1}(ii) do not include the assumption that $0\in  \intr{\mbox{dom}\, \A}.$ The reason for this difference is that the only role of the assumption $0\in  \intr{\doma}$ in Theorem~\ref{G-E} is to ensure that the upper bound~\eqref{upper_LDP} holds for {\em unbounded\/} sets~$B$ (cf.\ Corollary~6.1.6 in~\cite{DeZe09} or p.\,84 in~\cite{BoMo92}), whereas Theorem~\ref{t1}(ii) only asserts that the {\em local\/} LDP is satisfied for  $\{\ze_T\}$.

Note also that the smoothness condition in Theorem~\ref{t1}(ii) is only needed for obtaining the desired lower bound for  probabilities of the form $\pr (\ze_T\in (\al)_\ep ),$ as is seen from the theorem's proof in Section~\ref{Sec_2}.
\end{Rem}

\begin{Rem}
\rm
Remark~\ref{rem_A=A} shows that the claim of Theorem~\ref{t1}(ii) is consistent with the assertion of the classical G\"artner--Ellis Theorem~\ref{G-E}: if there exists an essentially smooth FF~$\Ao$ for $\{\ze_T\}$ then the WSFF~$A$ also exists and, in view of~\eqref{A=A}, both theorems assert the LLDP with a common rate function~$D ={\mathcal L}_\A ={\mathcal L}_\Ao .$
\end{Rem}

\begin{Rem}
\label{rem_SBF_3}
\rm
As   noted above, 
the WSFF  $\A$ is always automatically convex by virtue of its definition. However, the rate function~$D$ in  Theorem~\ref{t1}(i) does not need to be convex. The following examples illustrate the situation where the LLDP is satisfied with a non-convex rate function~$D$ and the WSFF~$\A={\mathcal L}_D$ exists (as it should, by Theorem~\ref{t1}(i)), but   $D\neq {\mathcal L}_\A$. It follows from  the general properties of the Legendre--Fenchel transform that, in such situations, ${\mathcal L}_\A$ is the largest convex lower semicontinuous minorant of~$D$ (see e.g.~\cite{Ro70}).

\end{Rem}

\begin{Example}
\label{Ex_2}
\rm
Suppose  that, for a function $\{a_T\}_{T\ge 0}$, one has
\[
\pr   (\ze_T=\al ) =\left\{
\begin{array}{ll}
2^{-T}, & \al =1,\\
1-2^{-T}, & \al =a_T,\\
0, & \al \not\in \{1,a_T\},
\end{array}
\right.
\quad T\ge 0.
\]

(i)~First assume that $a_T\to 0$ as $T\to\infty.$
One can  verify directly  that the LLDP is satisfied in this case, with the non-convex rate function
\[
D(\al)
=\left\{
\begin{array}{ll}
0, & \al =0,\\
\ln 2, & \al =1,\\
\infty, & \al \not\in \{0,1\},
\end{array}
\right.
\]
and that there exist both FF and WSFF equal to
\[
\Ao (\mu)= \A(\mu)=\max\{0,\mu - \ln 2\} = \mathcal L_D(\mu),\quad \mu \in \R.
\]
However, as $A$ is not smooth, Theorem~\ref{t1}(ii) is inapplicable. One can see from the proof in Section~\ref{Sec_2} that the ``jump'' of $A'(\mu)$ at $\mu =\ln 2$ makes it impossible to get the needed lower bound in terms of $\mathcal L_A$ for $\pr (\ze_T\in (\al)_\ep)$ when $\al \in (0,1).$ Hence we cannot claim that $D=\mathcal L_A$, and, indeed, an elementary calculation yields that
\[
{\mathcal L}_\A (\al)
 = \left\{
 \begin{array}{ll}
 \al \ln 2, & \al\in [0,1],\\
 \infty & \mbox{otherwise,}
 \end{array}
 \right.
\]
which is the largest convex lower semicontinuous minorant of~$D$, as it should be.

(ii)~Now assume that $|a_T|\to \infty$ as $T\to\infty.$ Clearly, the LLDP is also satisfied in this case, now with the convex rate function
\[
D(\al)
=\left\{
\begin{array}{ll}
\ln 2, & \al =1,\\
\infty, & \al \neq 1.
\end{array}
\right.
\]
The FF does not exist here, but the WSFF does and is equal to $\A(\mu)= \mu - \ln 2  = \mathcal L_D(\mu),$ $\mu \in \R,$ which satisfies the conditions of Theorem~\ref{t1}(ii).  Clearly, $\mathcal L_A(\al ) = D(\al),$ $\al \in \R.$

\end{Example}

The next example is slightly more sophisticated. It was obtained by modifying a construction from p.\,75 in~\cite{DeZe09}.

\begin{Example}\label{Ex_3}
\rm
Let $\{Y_t\}_{t\ge 0}$ be a discrete time Markov chain with the state space $\{1,2,3\}$,  transition probability matrix
  \[
  \left(
  \begin{array}{ccc}
  1/2 & 0 & 1/2 \\
  0 & 1/2 & 1/2 \\
  0& 0 & 1
  \end{array}
  \right),
  \]
and initial distribution $\pr (Y_0=1)=\pr (Y_0=2)=1/2.$ Set
\[
\xi_t:= (\ind (Y_t =1), \ind (Y_t =2))\in \R^2,
\quad
 \zeta_T  =(\zeta_{T,1}, \zeta_{T,2}) := \frac1T  \sum_{t=1}^T \xi_t ,
\quad
T=1,2,\ldots
\]
Clearly, $\pr (\zeta_{T,1}> 0, \zeta_{T,2}>0)=0$ for all $T> 0$  and, for any fixed $\al \in (0,1)$ and $\ep \in (0,\min\{\al, 1-\al\})$, one has
\begin{align*}
\pr (\ze_{T,j}\in (\al)_\ep )
  & = \sum_{k=\lfloor (\al-\ep) T\rfloor +2 }^{\lceil (\al+\ep) T\rceil  } 2^{-k} = 2^{-\lfloor (\al-\ep) T\rfloor +1} (1 +o(1)), \quad j=1,2.
\end{align*}
It is easily seen from here that, for the family   $\{\ze_T\}$, the LLDP is saisfied with the non-convex rate function
\begin{align*}
D(\al )
= \left\{
\begin{array}{ll}
  \al_1\ln  2,  & \al\in [0,1]\times\{0\},\\
\al_2\ln  2,  & \al\in \{0\}\times[0,1],\\
 \infty &  \mbox{otherwise}.
  \end{array}
\right.
\end{align*}
Direct computation yields
\begin{align}
\label{MC_L_D}
{\mathcal L}_D (\mu) = \max\{0, \max\{ \mu_1,\mu_2\} -\ln 2\},
\quad \mu =(\mu_1, \mu_2)\in \R^2.
\end{align}

On the other hand,
\begin{align*}
\exn e^{T \langle   \mu, \zeta_T\rangle}
& = \frac12\biggl(\sum_{k=0}^T e^{\mu_1 k }  2^{-k}
+ \sum_{k=0}^T    e^{ \mu_2 k}   2^{-k}\biggr),
\quad T= 1,2,\ldots
\end{align*}
Elementary calculation shows that the limit
$
\lim_{T\to\infty} \frac1T \ln \exn e^{T \langle \mu, \zeta _T\rangle}
$   exists and is equal to the right-hand side of~\eqref{MC_L_D}.
Therefore, there exist both the FF and WSFF (note that $|\zeta_T|\le 1$ always, meaning that these two functions are the same) given by
\[
\Ao (\mu ) = \A(\mu ) = {\mathcal L}_D (\mu),\quad \mu\in \R^d,
\]
in accordance with the claim of Theorem~\ref{t1}(i).

On the other hand, as $A$ is not smooth, one cannot claim that $D=\mathcal L_A.$  It is not hard to verify that
\[
{\mathcal L}_{\A} (\al)=
\left\{
\begin{array}{ll}
(\al_1 + \al_2) \ln 2, & \mbox{if  $\al_1+\al_2\le 1,$ $\al_1\ge0, $ $\al_2 \ge 0,$}\\
\infty & \mbox{otherwise}.
\end{array}
\right.
\]
Clearly, ${\mathcal L}_{\A} \neq D$ is  the largest convex lower semicontinuous minorant of~$D$, as it should be according to Remark~\ref{rem_SBF_3}.
\end{Example}

\section{Proof of Theorem~\ref{t1}}\label{Sec_2}

\Proof~of Theorem~\ref{t1}.
(i)~Let the LLDP with a rate function~$D$ be satisfied for $\{\zeta_T\}$. We will show that the WSFF exists for that family and that it is equal to    $\Aw :={\mathcal L}_D$.

First recall that the Legendre--Fenchel transform of any function is convex (see e.g.\ p.\,104 in~\cite{Ro70}), so that $\Aw$ is convex.

Next, assume that    $\mu \in  \mbox{dom} \,\Aw  $. In this case, for any $\eta >0$, there exists an $  \al_\eta (\mu)\in \R^d $ such that
\begin{align}
\label{>eta}
\langle \mu,  \al_\eta (\mu) \rangle  - D(\al_\eta (\mu)) >\Aw (\mu ) - \eta.
\end{align}
It is clear that, for any fixed $\ep>0$ and for $M\to \infty $ as $T\to\infty$, for all large enough~$T$ one has
$
(\al_\eta(\mu))_\ep \subset \Bom  .
$
Therefore, using~\eqref{lpbu}, we obtain, for  $\ep =\ep_T=\overline{o}(1)$, that
\begin{align*}
E_{T,M} (\mu)
 &  \ge \int_{(\al_\eta (\mu))_\ep} e^{T \langle \mu, x\rangle}  \pr (\ze_T \in dx)  
\notag \\
& \ge
\exp\{-T |\mu|  \ep + T  \langle \mu , \al_\eta (\mu) \rangle\}
 \pr \big(\ze_T \in (\al(\mu))_\ep\big)
\notag \\
& \ge  \exp\{- T |\mu|  \ep + T[\langle \mu , \al_\eta (\mu) \rangle - D(\al (\mu)) +o(1)] \}
\notag
\\
& \ge  e^{T(\Aw (\mu)-\eta  + o(1))}
\end{align*}
from~\eqref{>eta}.  Taking $\eta=\eta_T=\oo (1)$ and using the stnadard ``diagonalization argument'' (see e.g.\ \S\,7.16 in~\cite{RoFi10}), we get the lower bound
\begin{align}
E_{T,M} (\mu)
\ge
e^{T(\Aw (\mu)   + o(1))}.
\label{lower_b_BF}
\end{align}

Now assume that $\mu \not\in  \mbox{dom} \,\Aw  $, so that    $\Aw (\mu)=\infty$. Then, for any $N>0$, there will exist an $\al^N(\mu)$ such that
\[
\langle \mu,  \al^N (\mu) \rangle  - D(\al^N (\mu)) > N.
\]
Arguing as above, we will show that $E_{T,M} (\mu)\ge e^{T(N + o(1))}$. Since $N$ is arbitrary, this means that  $T^{-1}\ln E_{T,M} (\mu)\to \infty $  as  $T\to \infty$. Therefore we established that~\eqref{2} holds true in this case.

It  is clear that we only need to obtain  an upper bound for $E_{T,M} (\mu) $ in the case when $\Aw (\mu ) <\infty.$  Consider first a fixed  $M>0$. By virtue of~\eqref{loc_ld_sup}, for any $\delta >0$ and $\al \in \Bom$,
there exists an $\ep (\al, \delta) \in (0,\delta)$  such that
\[
\limsup_{T\to\infty} \frac1T \ln \pr (\ze_T \in (\al )_{\ep (\al, \delta)})
 \le -D_\delta (\al )+\delta,
\]
where we used $D_\delta (\al ) := \min\{D(\al), \delta^{-1}\}$ rather than just $D(\al)$  to ensure that this inequality  will hold in the case when $D(\al)=\infty$ as well.
For a fixed  $\delta >0,$ the neighbour\-hoods   $(\al)_{\ep (\al, \delta)}$, $\al \in \Bom ,$ form an open cover of the compact ball  $\Bom $. Hence there exists a finite subcover   $\{(\al)_{\ep (\al^{(k)} , \delta)}\}_{k\le n_\delta}$ thereof. For all large enough~$T$ one has
\begin{align*}
\max_{k\le n_\delta } \big[ T^{-1}\ln \pr \big(\ze_T \in (\al^{(k)})_{\ep (\al^{(k)}, \delta)}\big)+ D_\delta (\al^{(k)})\big] \le  2\delta,
\end{align*}
which implies that
\begin{align*}
E_{T,M} (\mu)
 &  \le \sum_{k=1}^{n_\delta}  \int_{(\al^{(k)})_{\ep (\al^{(k)}, \delta)}} e^{ T \langle \mu ,   x\rangle } \pr (\ze_T\in dx)
 \\
 &
 \le  \sum_{k=1}^{n_\delta} e^{T \langle \mu,  \al^{(k)}\rangle + T |\mu|  \ep (\al^{(k)}, \delta)}  e^{- T (D_\delta(\al^{(k)})-2\delta)}
\\
& \le  e^{(|\mu|+2)T\delta } n_\delta \exp\big\{\max_{k\le n_\delta}(\langle \mu,\al^{(k)}\rangle- D_\delta (\al^{(k)}))\big\}
\\
& =  e^{(|\mu|+2)T\delta } n_\delta \exp\big\{\max_{k\le n_\delta}(\langle \mu,\al^{(k)}\rangle- D  (\al^{(k)}))\big\}
\\
&   \le  e^{(|\mu|+2)T\delta } n_\delta e^{T \Aw (\mu)  },
\end{align*}
where the validity of the second last relation (the equality) is guaranteed under the assumption that  $\delta^{-1}\ge M|\mu|- \Aw (\mu),$ which does not restrict the generality.

It remains to choose  $M=M_T\to\infty$ and $\delta = \delta_T\to 0$ so slowly as $T\to\infty$ that the above bound will still hold true and one would have  $n_\de =e^{o(T)}.$ This, together with~\eqref{lower_b_BF}, establishes the existence of the WSFF $\A=\Aw $.

\medskip

{\rm (ii)}~Assume now that there exists an essentially smooth WSFF  $\A(\mu)$.  
We will prove that the LLDP is satisfied for $\{\zeta_T\}$ with the rate function equal to $\Dw :=\mathcal L_{\A}$.

{\em Upper bound.} Take an arbitrary $\al\in \R^d$ and let $M,\ep>0$ be such that  $(\al)_\ep\subset \Bom $. Then, for  any $\mu \in \R^d$, one has
\begin{align}
\label{ind_ind}
\ind(\ze_T\in (\al)_\ep )\le e^{ T \langle \mu, (\ze_T -\al)\rangle + T |\mu| \ep }\ind(\ze_T\in \Bom  ).
\end{align}
Hence, for  $M=M_T \to\ui$ and  $\ep=\ep_T =\oo (1)$ as  $T\to\infty$, we obtain that
\begin{align*}
\pr (\ze_T \in (\al)_\ep )
& \le e^{- T \langle \mu, \al\rangle + T|\mu|\ep  } \exn (e^{ T\langle\mu, \ze_T\rangle}; \ze_T \in  \Bom  )
\\
& = \exp\{ -T(\langle\mu ,\al \rangle- A(\mu)+o(1))\}.
\end{align*}
Since  $\mu $ is arbitrary, we thereby showed that
\begin{align*}
\pr (\ze_T \in (\al)_\ep )
 \le   \exp\{ -T( \widehat D(\al) +o(1))\},
\end{align*}
which immediately implies~\eqref{loc_ld_sup} with $D=\Dw$.

{\em Lower bound.}
Here  we will derive a lower bound for the probability on the left-hand side of the above formula (note that we only need it for $\al\in \mbox{dom}\, \Dw$). To this end, we will employ an idea from the proof of Theorem~2.3.6 in~\cite{DeZe09}. Namely, we will apply the already obtained  upper bound after a suitable  change of measure.

Denote by~$P_T$ the distribution of $\ze_T$ and, for given $M>0$ and $\mu\in \R^d$, consider a new distribution~$\widetilde{P}_{T,M,\mu}$ on~$\R^d$, specifying it via its density
\begin{align}
\label{density}
\frac{d\widetilde{P}_{T,M,\mu}}{dP_T}(x)
 := \frac{e^{T \langle \mu,  x\rangle }}{E_{T,M }(\mu)}\ind (x\in \Bom ).
\end{align}
Then, for   $\al\in \R^d$ and $\ep>0$ such that $(\al)_\ep \subset  \Bom $, one has
\[
P_T((\al)_\ep) =E_{T,M}(\mu)\int_{(\al)_\ep}e^{-T \langle \mu ,  x\rangle }d\widetilde{P}_{T,M,\mu}(x),
\]
so that
\begin{align}
\frac1T \ln  P_T((\al)_\ep)
&
= \frac1T  \ln E_{T,M}(\mu) - \langle \mu , \al \rangle
 + \frac1T\ln \int_{(\al)_\ep}e^{ T\langle \mu,   \al - x \rangle }d\widetilde{P}_{T,M,\mu}(x)
\notag \\
& \ge
 \frac1T  \ln E_{T,M}(\mu) - \langle\mu, \al \rangle-|\mu| \ep
 + \frac1T\ln  \widetilde{P}_{T,M,\mu}( (\al)_\ep).
 \label{repa}
\end{align}
Therefore, for any fixed  $\al\in \R^d$ and for $M=M_T\to\ui$, $\ep=\ept =o(1),$ one has from~\eqref{2} that
\begin{align}
\lim_{\ep \searrow  0}\liminf_{T\to \infty} \frac1T & \ln  P_T((\al)_\ep)
\ge  \A (\mu)  - \langle\mu, \al\rangle
+ \lim_{\ep \searrow 0}\liminf_{T\to \infty}  \frac1T\ln  \widetilde{P}_{T,M,\mu}( (\al)_\ep)
\notag
\\
&
\ge -   \Dw (\al) +  \lim_{\ep \searrow  0}\liminf_{T\to \infty}  \frac1T\ln [1- \widetilde{P}_{T,M,\mu}( (\al)_\ep^c)],
\label{lower_bound_1}
\end{align}
where the last inequality holds since $\mu\in\R^d$ is arbitrary.

To establish the desired lower bound, it remains to show that, for a suitably chosen $\mu = \mu (\al),$
\[
\widetilde{P}_{T,M,\mu}( (\al)_\ep^c)=o(1)\quad \mbox{ for  $\ep=\ep_T= \overline{o}(1).$}
\]

Due to the assumption on the existence of the WSFF, for any fixed  $\la \in \R^d$ and for $M=M_T\to \ui$, one has
\begin{align}
\frac1T \ln E_{T,M}(\mu+\la) = A(\mu +\la) + \theta_{\mu+\la} (T),
\label{E_A}
\end{align}
where $\theta_{\mu+\la} (T) \to 0$ as $T\to\infty$, so that
\begin{align}
\int e^{\langle \la, T x\rangle} d\widetilde{P}_{T,M,\mu}(x)
 &
 =  \frac1{{E_{T,M}(\mu)}} \int_{\Bom } e^{ T \langle (\mu +\la),   x\rangle}d {P}_{T }(x)
 \notag
 \\
 &
 = \frac{E_{T,M}(\mu+\la)}{E_{T,M}(\mu)}
 \notag
 \\
 &
 = \exp\{T(\A(\mu+\la) - \A(\mu) + \theta_{\mu+\la} (T) -\theta_{\mu}  (T))\}.
 \label{lambda_mom}
\end{align}

Let  $\widetilde{\ze}_{T } =(\widetilde{\ze}_{T,1} , \ldots, \widetilde{\ze}_{T,d} )$ be a random vector with distribution $\widetilde{P}_{T,M,\mu}$. For $\al = (\al_1, \ldots, \al_d)$  and $\ep_1:= \ep d^{-1/2}$, one clearly has
\begin{align}
\label{bound_sum}
\widetilde{P}_{T,M,\mu}( (\al)_\ep^c)
 = \pr (\widetilde{\ze}_{T } \not\in (\al)_\ep)
 \le \sum_{j=1}^d \pr \big(|\widetilde{\ze}_{T,j} -\al_j|\ge \ep_1 \big).
\end{align}
Setting $\la (\de):=(\de, 0, \ldots,0)\in \R^d$ for $\de >0,$ from  the exponential Chebyshev inequality and~\eqref{lambda_mom}  we get
\begin{align}
\pr (\widetilde{\ze}_{T,1} & \ge \al_1+\ep_1)
  \le \frac{\exn e^{T \de \widetilde{\ze}_{T,1}  }}{e^{ T \de    (\al_1 +\ep_1)}}
  =\frac{\exn e^{ T\langle \la  (\de),  \widetilde{\ze}_{T} \rangle  }}{e^{T \de   (\al_1 +\ep_1)}}
  \notag\\
  &
 = \exp\{T(\A(\mu+\la (\de)) - \A(\mu) - \de   (\al_1 +\ep_1) + \theta_{\mu+\la(\de)} (T) -\theta_{\mu}  (T))\}.
 \label{upper_exp}
\end{align}

First assume that $\al \in \intr{{\rm dom}\, \Dw }$. Since $\A$ is essentially smooth, by Theorem~26.3 in~\cite{Ro70}  there exists a $\mu =\mu (\al)  \in  \intr{{\rm dom}\, \A }$  such that
\begin{align}
\grad \A(\nu) |_{\nu = \mu(\al) } =\al.
 \label{grad_A}
\end{align}
The function $a(h):= \A(\mu+\la(h)),$ $h\in \R,$ is clearly convex. Since $a'(0) =\al_1$ by virtue of~\eqref{grad_A}, we conclude that
\[
\A(\mu+\la(\de ) ) - \A(\mu) =  a(\de) - a(0) = \al_1\de +b(\de)\de , \quad \de >0,
\]
where $ \widehat b (h):=\sup_{0\le \de \le h}| b (\de)|=o(1)$ as $h \to 0$.
Now~\eqref{upper_exp} implies that
\begin{align}
\pr (\widetilde{\ze}_{T,1} \ge \al_1+\ep_1)
 &
  \le   \exp\{- T(   \de   \ep_1- b(\de) \de -  \theta_{\mu+\la(\de)} (T) +\theta_{\mu}  (T))\}.
 \label{upper_exp_1}
\end{align}
We will show  that one can choose  $\ep,  \de=\oo (1)$ as $T\to \infty$ such that the argument in the above exponent will tend to~$-\infty$.

Take a  $\de^{(0)}>0$ such that $\mu + \la(\de ^{(0)})\in \intr{{\rm dom}\, A}$  and put  $\de^{(k)}:= \de^{(0)}/k$,  $k=1,2,\ldots$ Further, one can always find  $0=T_0<T_1 <T_2 <\cdots$ such that
\[
\sup_{T\ge T_k} (| \theta_{\mu+\la(\de^{(k)})} (T)| + |\theta_{\mu}  (T)|)
 \le   \widehat b (\de^{(k)})\de^{(k)} , \quad k\ge 1.
\]
Setting
\begin{align*}
\de
& =\de_T:= \sum_{k\ge 0}  \de^{(k)} \ind (T_k \le T<T_{k+1}),
\\
\ep
& =\ep_T:= \sum_{k\ge 0}    \widehat b (\de^{(k)})^{1/2} \ind (T_k \le T<T_{k+1}),
\end{align*}
we obtain the desired construction, thus showing that
\[
\pr (\widetilde{\ze}_{T,1} \ge \al_1+\ep_1)\to 0.
\]
That $\pr (\widetilde{\ze}_{T,1} \le \al_1-\ep_1)\to 0$ for $\ep=\ept= \oo (1) $ is shown in the same way, using $\de <0.$

Since the same argument works for all the other components $\widetilde{\ze}_{T,j}$, $1<j\le d,$ we see from~\eqref{bound_sum} that  $\widetilde{P}_{T,M,\mu}( (\al)_\ep^c)\to 0.$ This, together with~\eqref{lower_bound_1}, establishes the desired lower bound in the case when $\al \in \intr{{\rm dom}\, \Dw }$.

Finally, we consider the alternative situation where $\al \in ( \partial {\rm dom}\, \Dw )\cap ( {\rm dom}\, \Dw )$. For simplicity, we will only consider the case $d=1$ (the general case is dealt with in the same way, the argument becoming somewhat more cumbersome).

One can assume without loss of generality that $\al = \max ( {\rm dom}\, \Dw )$.
Take an $\ep >0$ such
\[
\al':= \al - \ep' \in \intr{{\rm dom}\, \Dw} , \quad \mbox{where}\  \ep' := \ep/2.
\]
Using~\eqref{repa} and~\eqref{E_A}, it is easily seen that
\begin{align*}
\frac1T \ln  P_T ((\al)_\ep)
&
\ge \frac1T \ln  P_T((\al')_{\ep'})
 \\
& \ge
 \frac1T  \ln E_{T ,M}(\mu) - \mu \al'   -|\mu| \ep'
 + \frac1T\ln  \widetilde{P}_{T, M, \mu} ( (\al')_{\ep'})
 \\
 & \ge
 - \Dw (\al) - |\mu| \ep + \theta_\mu (T)
  + \frac1T\ln  \widetilde{P}_{T, M, \mu} ( (\al')_{\ep'}).
\end{align*}
Since  $\al' \in \intr{{\rm dom}\, \Dw} ,$  the argument we used above in the case when  $\al\in  \intr{\domd} $ (now with the choice  $\mu = \mu (\al')$)   shows that, for small enough  $\de$, one has
\begin{align*}
\frac1T \ln  P_T((\al)_\ep)
\ge
 -  \Dw (\al) - 2 |\mu| \ep
 \end{align*}
for all large enough~$T$. Since~$\ep >0$ is arbitrary, this implies the validity of the desired lower bound. Theorem~\ref{t1} is proved.
\hfill$\Box$

\bigskip\noindent
{\bf Acknowledgments.} The author is grateful to A.~A.~Borovkov for drawing his attention to the problem and for valuable discussions.


\begin{thebibliography}{99}

\bibitem{Bo20a}
{\it Borovkov, A.\,A.}
{Asymptotic Analysis of Random Walks: Light-Tailed Distributions.}
{Cambridge: Cambridge Univ.\ Press, 2020.}


\bibitem{Bo20} 
{\it Borovkov, A.\,A.}
{Compound Renewal Processes.}
Cambridge Univ.\ Press
2022. 




\bibitem{BoMo92}
{\it Borovkov, A.\,A., Mogulskii, A.\,A.}
{Large deviations and testing of statistical hypotheses. I. Large deviations of sums of random vectors.}
{\em Siberian Adv. Math.} 1992, {\bf 2}, 52--120.

\bibitem{BoMo10}
{\it Borovkov, A.\,A., Mogulskii, A.\,A.}
On large deviation principles in metric spaces.
{\em Siberian Math. J.}  2010, {\bf 51}, 989--1003.




\bibitem{BoMo11}
{\it Borovkov, A.\,A., Mogulskii, A.\,A.}
{On large deviation principles for random walk trajectories. I.}
{\em Theor.\ Probab.\ Appl.} 2011, {\bf 56}, 538--561.



\bibitem{Cr38}
{\it Cram\'er, H.}
Sur un nouveau th\'eor\`eme–limite de la th\'eorie des probabilit\'es. In: {\em Actualit\'es Scientifiques et Industrielles}, number 736 in Colloque consacr\'e \`a la th\'eorie des probabilit\'es, pp. 5--23. Hermann, Paris, 1938. [English translation: {\tt arXiv:1802.05988v4}.]


\bibitem{DeZe09}
{\it Dembo\,A., Zeitouni. O.}
Large Deviations Techniques and Applications. Corr.\ printing of 2nd edn., Springer, 2009.

\bibitem{DeSt89}
{\it Deuschel, J.\,D., Stroock, D.\,W.}
{Large Deviations.} Academic Press, Boston, 1989.


\bibitem{El84}
{\em Ellis, R. S.}
Large deviations for a general class of random vectors.
{\em Ann.\ Probab.} 1984, 12, 1--12.

\bibitem{Ga77}
{\em G\"artner, J.}
On large deviations from the invariant measure.
{\em Theor.\ Prob.\ Appl.} 1977, 22, 24--39.

\bibitem{Pu01}
{\em Puhalskii, A.}
{Large Deviations and Idempotent Probability.}
{Boca Raton, FL: CRC Press, 2001.}


\bibitem{Ro70}
{\em Rockafellar, R.T.}
{Convex Analysis.} {Princeton: Princeton Univ.\ Press, 1970. }


\bibitem{RoFi10}
{\em Royden~H.L., Fitzpatrick~P.M.}
{Real Analysis.} 4th ed.
{Boston: Pearson,  2010.}


\bibitem{Va66}
{\em Varadhan, S. R. S.}
Asymptotic probabilities and differential equations.
{\em Comm. Pure Appl. Math.} 1966, 19, 261--286.
\end{thebibliography}
\end{document}